\documentclass[10pt,twoside]{article}
\usepackage{graphicx,amsmath,amsfonts}
\usepackage{Latex-document}

\def\ll{\:\raisebox{-0.6ex}{$\stackrel{=}{def}$}\:}
\def\CC{{\mathbb C}}
\def\PP{{\mathbb P}}
\def\ZZ{{\mathbb Z}}

\title
{\bf Representations of Algebraic Groups and\vskip -2mm Principal
Bundles on Algebraic Varieties\vskip 6mm}

\author{Vikram Bhagvandas Mehta\vspace*{-0.5cm}\thanks{Tata Institute of Fundamental Research, Mumbai, India. E-mail:
vikram@math.tifr.res.in}}

\date{\vspace{-8mm}}

\begin{document}

\maketitle

\thispagestyle{first} \setcounter{page}{629}

\begin{abstract}

\vskip 3mm

In this talk we discuss the relations between representations of
algebraic groups and principal bundles on algebraic varieties,
especially in characteristic $p$. We quickly review the notions of
stable and semistable vector bundles and principal $G$-bundles ,
where $G$ is any semisimple group. We define the notion of a low
height representation in characteristic $p$ and outline a proof of
the theorem that a bundle induced from a semistable bundle by a
low height representation is again semistable. We include
applications of this result to the following questions in
characteristic $p$:

1) Existence of the moduli spaces of semistable $G$-bundles on
curves.

2) Rationality of the canonical parabolic for nonsemistable
principal bundles on curves.

3) Luna's etale slice theorem.

We outline an application of a recent result of Hashimoto to study
the singularities of the moduli spaces in (1) above, as well as
when these spaces specialize correctly from characteristic $0$ to
characteristic $p$. We also discuss the results of
Laszlo-Beauville-Sorger and Kumar-Narasimhan on the Picard group
of these spaces. This is combined with the work of Hara and
Srinivas-Mehta to show that these moduli spaces are $F$-split for
$p$ very large. We conclude by listing some open problems, in
particular the problem of refining the bounds on the primes
involved.

\vskip 4.5mm


\end{abstract}

\vskip 12mm



\section{Some Definitions}

\vskip-5mm \hspace{5mm}

We begin with some basic definitions:\\

Let $V$ be a vector bundle on a smooth projective curve $X$ of
genus $g$ over an algebraically closed field (in any
characteristic).

\noindent {\bf Definition 1.1:} {\it $V$ is \emph{stable} (\emph{
respectively semi-stable}) if for all subbundles $W$ of $V$, we
have
$$\mu(W) \ll \emph{deg} \ W/rk \ W \  < (\leq)$$
$$\mu(V) \ll \ \emph{deg} \ V/rk \ V.$$
}

\noindent For integers $r$ and $d$ with $r > 0$, one constructs
the moduli spaces $U^s (r,d) (U(r,d))$ of stable (semistable)
vector bundles of rank $r$ and degree $d$, using Geometric
Invariant Theory (G.I.T.).

If the ground field is $\CC$, the complex numbers, one has the
basic (genus $X \ge 2$):

\noindent {\bf Theorem 1.2:} {\it Let $V$  have degree $0$. Then
$V$ is stable $\Leftrightarrow V \simeq V_{\sigma}$, for some
irreducible representation $\sigma : \pi_1(X) \to U(n)$.}

This is due to Narasimhan-Seshadri.  Note that $H \to X$ is a
principal $\pi_1(X)$ fibration, where $H$ is the upper-half plane.
Any $\sigma : \pi_1(X) \to GL(n, \CC)$ gives a vector bundle of
rank $n$ on $X, V_\sigma = H \times^{\pi_1(X)} \CC^n$.

\noindent {\bf Remark 1.3:}  It follows from Theorem 1.2 that if
$V$ is a semistable bundle on a curve $X$ over $\CC$, then
$\otimes^n(V), S^n(V)$, in fact any bundle induced from $V$ is
again semistable.  By Lefschetz, this holds for any algebraically
closed field of characteristic $0$.

\noindent {\bf Remark 1.4:} In general, a subbundle $W$ of a
vector bundle $V$ is a reduction of the structure group of the
principal bundle of $V$ to a maximal parabolic of $GL(n), n =
\emph{rank} \ V$.  This is in turn equivalent to a section
$\sigma$ of the associated fibre space:
$$E \times^{GL(n)} \ \ GL(n)/P. $$
Now let $X$ be a smooth curve and $E\stackrel{\pi}{\rightarrow} X$
a principal $G$-bundle on $X$, where $G$ is a semisimple (or even
a reductive) group in any characteristic.

\noindent {\bf Definition 1.5:} $E$ is stable (semistable)
$\Leftrightarrow \forall$ maximal parabolics $P$ of $G, \forall$
sections
 $\sigma$ of
$E(G/P)$, we have degree $\sigma^\# T_{\pi} > 0 (\ge 0)$, where
$T_{\pi}$ is the relative tangent bundle of $E(G/P)
\stackrel{\pi}{\rightarrow} X$.

Over $\CC$, we have the following [18]:

\noindent {\bf Theorem 1.6:} {\it $E \to X$ is stable
$\Leftrightarrow E \simeq E_\sigma$ for some irreducible
representation $\sigma : \pi_1(X) \rightarrow K$, the maximal
compact of $G$}.

The analogue of Remark 1.3 is valid in this general situation.

\noindent {\bf Remark 1.7:}  One can analogously define stable and
semistable vector bundles and principal bundles  on normal
projective varieties of dimension $>1$.  Again, in characteristic
$0$, bundles induced from semistable bundles continue to be
semistable.

\noindent {\bf Remark 1.8:}  In characteristic $p$, bundles
induced from semistable bundles need not be semistable, in
general[7]. In this lecture we shall examine some conditions when
this does hold, and also discuss some applications to the moduli
spaces of principal $G$-bundles on curves.


\section{Low height representations}

\vskip-5mm \hspace{5mm}

Here we introduce the basic notion of a low height representation
in characteristic $p$.  Let $f:G \to SL(n) = SL(V)$ be a
representation of $G$ in char $p$, $G$ being reductive.  Fix a
Borel $B$ and a Torus $T$ in $G$. Let $L(\lambda_i), 1 \leq i \leq
m$, be  the simple $G$-modules occurring in the Jordan-Holder
filtration of $V$. Write each $\lambda_i$ as
$\displaystyle{\sum_{j}} q_{ij} \alpha_j$, where $\{ \alpha_j \}$
is the system of simple roots corresponding to $B$ and $q_{ij} \in
Q \ \forall i,j$.  Define $ht \lambda_i = \displaystyle{\sum_j}
q_{ij}$. Then one has the basic [9,20]:

\noindent {\bf Definition 2.1:}  {\it $f$ is a low-height
representation of $G,$ or $V$ is a low-height module over $G$, if
$2 ht (\lambda_i) < p \  \forall i$.}

\noindent {\bf Remark 2.2:}  If $2 ht (\lambda_i)< p \ \forall i$,
then it easily follows that $V$ is a completely reducible
$G$-module.  In fact for any subgroup $\Gamma$ of $G, V$ is
completely reducible over $\Gamma \Leftrightarrow \Gamma$ itself
is completely reducible in $G$. By definition, an abstract
subgroup $\Gamma$ of $G$ is completely reducible in $G
\Leftrightarrow$ for any parabolic $P$ of $G$, if $\Gamma$ is
contained in $P$ then $\Gamma$ is contained in a Levi component
$L$ of $P$. These results were proved by Serre[20] using the
notion of a saturated subgroup of $G$.

In general, denote sup $(2 ht \ \lambda_i)$ by $ht_GV$.  If $V$ is
the standard $SL(n)$ module, then $ht_{SL(n)} \wedge^i (V)= i
(n-i), 1 \leq i \leq n-1$. More generally, $ht_G (V_1 \otimes V_2)
= ht_G V_1+ht_G V_2$.  The following theorem is the key link
between low-height representations and semistability of induced
bundles [9]:

\noindent {\bf Theorem 2.3:}  \emph{Let $E \to X$ be a semistable
$G$-bundle, where $G$ is semisimple and the base $X$ is a normal
projective variety. Let $f:G \to SL(n)$ be a low-height
representation.  Then the induced bundle $E (SL(n))$ is again
semistable}.

The proof is an interplay between the results of Bogomolov, Kempf,
Rousseau and Kirwan in G.I.T. on one hand and the results of Serre
mentioned earlier on the other.  The group scheme $E(G)$ over $X$
acts on  $E(SL(n)/P)$ and assume that $\sigma$ is a section of the
latter.  Consider the generic point $K$ of $X$ and its algebraic
closure $\overline{K}$.  Then $E(G)_{\overline{K}}$ acts on
$E(SL(n)/P)_{\overline{K}}$, and $\sigma$ is a $K$-rational point
of the latter.  There are 2 possibilities:

\begin{enumerate}
\item[1)]  $\sigma$ is $G.I.T$ semistable.  In this case, one can easily
prove that deg $\sigma^{\#} T_{ \pi} \geq 0$.

\item[2)]  $\sigma$ is $G.I.T.$ unstable, i.e., not semistable.  Let
$P (\sigma)$ be the Kempf-Rousseau parabolic for $\sigma$, which
is defined over $\overline{K}$.  For deg $\sigma^{\#} T_{\pi}$ to
be $\geq 0$ it is {\it sufficient} that $P (\sigma)$ is defined
over $K$. Note that since $V$ is a low-height representation of
$G$, one has $p \geq h$.  One then has ([20]).
\end{enumerate}

\noindent {\bf Proposition 2.4:}  \emph{ If $p \geq h$, there is a
unique G-invariant isomorphism log: $G^u \to \underline{g}_{{\rm
nilp}}$, where $G^u$ is the unipotent variety of $G$ and
$\underline{g}_{{\rm nilp}}$ is the nilpotent variety of
$\underline{g} =$ Lie $G$}.

Proposition 2.4 is used in

\noindent {\bf Proposition 2.5:}   \emph{Let $H$ be any semisimple
group and $W$ a low-height representation of $H$.  Let $ W_1
\subset W$ and assume that $\exists X \in $ Lie $H, \  X$
nilpotent such that $X \in $ Lie (Stab $(W_1))$.  Then in fact one
has $X \in $ Lie [Stab $(W_1)_{{\rm red}} ]$}.

Along with some facts from $G.I.T$, Proposition 2.5 enables us to
prove that $P (\sigma)$ is in fact defined over $K$, thus
finishing the sketch of the proof of Theorem 2.3.  See also
Ramanathan-Ramanan [19]. One application of low-height
representations is in the proof of a conjecture of Behrend on the
rationality of the canonical parabolic or the instability
parabolic.  If $V$ is a nonsemistable bundle on a variety $X$,
then one can show that there exists a flag $V^\cdot$,
$$0=V_0 \subset V_1 \subset V_2 \cdots \subset V_n =V$$
of subbundles of $V$ with the properties:

\begin{enumerate}
\item[(1)] Each $V_i / V_{i-1}$ is semistable and $\mu \ (V_i/V_{i-1})>
\mu \ (V_{i+1} /V_i), 1 \leq i \leq n-1$.
\item[(2)] The flag $V^\cdot$ as in (1) is unique and infinitesimally
unique, i.e., $V^\cdot$ is defined over any field of definition of
$X$ and $V$.  Such a flag corresponds to a reduction to a
parabolic $P$ of $GL(n)$ and properties (1) and (2) may be
expressed as follows: the {\it elementary} vector bundles on $X$
associated to $P$ all have positive degree and $H^0(X,
E(\underline{g}) / E(\underline{p}))=0$, where $\underline{g} =
\mbox{Lie} \ GL(n)$ and $\underline{p} = \mbox{Lie} \ P$.
\end{enumerate}

One may ask whether there is a such a canonical reduction for a
nonsemistable principal $G$ bundle $E \rightarrow X$.  Such a
reduction was first asserted first by Ramanathan [18], and then by
Atiyah-Bott[1] ,both over $\CC$ and both without proofs. It was
Behrend [ 5 ], who first proved the existence and uniqueness of
the canonical reduction to the instability parabolic in all
characteristics. Further, Behrend
 conjectured  that $H^0(X, E(\underline{g})/ E(\underline{p}))=0$.

In characteristic zero, one can check that all three definitions
of the instability parabolic coincide and that Behrend's
conjecture is valid.  In characteristic $p$, one uses low-height
representations to show the equality of the three definitions and
prove Behrend's conjecture [14].

\noindent {\bf Theorem 2.6:}  \emph{Let $E \rightarrow X$ be a
nonsemistable principal $G$-bundle in char $p$.  Assume that $ p >
2  dimG $. Then all the 3 definitions coincide and further we have
$H^0(X, E(\underline{g})/E(\underline{p})) =0$,  where
$\underline{p} =
\mbox{Lie} \ P$ and $P$ is the instability parabolic}.\\

Theorem 2.6 is useful, among other things, for classifying principal $G$-bundles on $\PP^1$ and $\PP^2$ in characteristic $p$.\\

If $V$ is a finite-dimensional representation of a semisimple
group $G$ (in any characteristic), then the G.I.T. quotient $V//G$
parametrizes the closed orbits in $V$. Now, let the characteristic
be zero and let $v_0 \in V$ have a closed orbit. Then Luna's
\'etale slice theorem says that $\exists$ a locally closed
non-singular subvariety $S$ of $V$ such that $v_0 \in S$ and $S //
G_{v_o}$ is isomorphic to $V//G$, locally at $v_0$, in the \'etale
topology.  Here $G_{v_0}$ is the stabilizer of $v_0$.  The proof
uses the fact that $G_{v_0}$ is a reductive subgroup of $G$ (not
necessarily connected!), hence $V$ is a completely reducible $G$
module.  In characteristic $p$, one has to assume that $V$ is a
low-ht representation of $G$.  Then the conclusion of Luna's
\'etale slice theorem is still valid: to be more precise, let $V$
be a low-ht representation of $G$ and let $v_0 \in V$ have a
closed orbit.  Put $H$ = Stab $(v_0)$.  The essential point, as in
characteristic 0, is to prove the complete reducibility of $V$
over $H$.  Using the low-ht assumption, one shows that every $X
\in \ \mbox{Lie} \ H$ with $X$ nilpotent can be integrated to a
homomorphism $G_a \rightarrow H$ with tangent vector $X$.  Now,
under the hypothesis of low-ht, one shows that $H_{\mbox{red}}$ is
a saturated subgroup of $G$ and $(H_{\mbox{red}} :
H^0_{\mbox{red}})$ is prime to  $p$. This shows that $V$ is a
completely reducible $H_{\mbox{red}}$ module.  Further, one shows
that $H_{\mbox{red}}$ is a normal subgroup of $H$  with
$H/H_{\mbox{red}}$ a finite group of multiplicative type, i.e. a
finite subgroup of a torus.  Now the complete reducibility of $V$
over $H$ follows easily [11]. Just as in characteristic zero,
 one deduces the existence of a smooth $H$-invariant subvariety $S$ of
$V$ such that $v_0 \in S$ and $S // H $ is locally isomorphic to
$V //G$ at $v_0$.  This result is used in the construction of the
moduli space $M_G$ to be described in the next section.

\section{Construction of the moduli spaces}

\vskip-5mm \hspace{5mm}

The moduli spaces of semistable $G$-bundles on curves were first
constructed by Ramanathan over $\CC$ [16,17], then by Faltings and
Balaji-Seshadri in characteristic 0 [3,6].  There are 3 main
points in Ramanathan's construction:

\begin{enumerate}
\item If $E \rightarrow X$ is semistable, then the adjoint bundle
$E(\underline{g})$ is semistable.
\item If $E \rightarrow X$ is polystable, then $E(\underline{g})$ is also
polystable.
\item A semisimple Lie Algebra in char 0 is rigid.
\end{enumerate}

The construction of $M_G$ in char $p$ was carried out in [2,15].
We describe the method of [15] first : points (1) and (2) are
handled by Theorem 2.3 and the following [11] :

\noindent {\bf Theorem 3.1:} \emph{Let $E \rightarrow X$ be a
polystable $G$-bundle
 over a curve in char $p$.  Let $\sigma : G \rightarrow SL(n) = SL(V) $ be a
 representation such that all the exterior powers
 $\wedge^i V, 1 \leq i \leq n-1$,  are low-height
representations.  Then the induced bundle $E(V)$ is also
polystable}.

The proof uses Luna's \'etale slice theorem in char $p$ and
Theorem 2.3.

Now one takes a total family $T$ of semistable $G$ bundle on $X$
and takes the
 good quotient of $T$ to obtain $M_G$ in char $p$.  Theorem 3.1 is used to
 identify the closed points of $M_G$ as the isomorphism classes of polystable
 $G$-bundles, just as in
char 0. The semistable reduction theorem is proved by lifting to
characteristic 0 and then applying Ramanathan's proof (in which
(3)
 above plays a crucial role).  This construction follows Ramanathan very
 closely and, as is clear, one has to make
low-height assumptions as in Theorem 3.1.

The method of [2] follows the one in [3] with some technical and
conceptual changes. One chooses an embedding $G \rightarrow SL(n)$
and a
 representation $W$ for $SL(n)$ such that (1) $G$ is the stabilizer of
 some $w_0 \in W $. (2) $W$ is a ``low separable index representation'' of
 $SL(n)$, i.e., all stabilizers are
reduced and $W$ is low-height over $SL(n)$.  The semistable
reduction theorem is proved using the theory of Bruhat-Tits.  Here
also suitable low-height assumptions have to be made.


\section{Singularities and specialization of the moduli spaces}

\vskip-5mm \hspace{5mm}

We first discuss the singularities of $M_G$, assuming throughout
that $G$ is simply connected.  In char 0, $M_G$ has rational
singularities, this follows from Boutot's theorem.  In char $p$,
the following theorem due to Hashimoto [8] is relevant:

\noindent {\bf Theorem 4.1:} \emph{Let $V$ be a representation of
$G$ such that all the symmetric powers $S^n(V)$ have a good
filtration.  Then the ring of invariant $[S^\cdot(V)]^G$ is
strongly $F$-regular}.

Strongly $F$-regularity is a  notion in the theory of tight
closure in commutative algebra.  We just note that if a geometric
domain is strongly $F$-regular then it is normal,Cohen-Macaulay,
$F$-split and has ``rational-like'' singularities. Now let $ t \in
M_G$ be the ``worst point'', i.e., the trivial $G$-bundle on $X$.

The local ring $({\mathcal O}_{M_G}, t)^\wedge$ is isomorphic to
$(S^\cdot(W)// G)^\wedge$, where $W$ = direct sum of $g$ copies of
$\underline{g}$, with $G$ acting diagonally. If $p$ is a good
prime for  $G$ , then Hashimoto's theorem implies that ${\mathcal
O}_{M_G,t}$ is strongly $F$-regular. The other points of $M_G$ are
not so well understood. This would require a detailed study of the
automorphism groups of polystable bundles, both in char 0 and $p$,
and of their invariants of the slice representations. This is
necessary also to study the specialization problem, i.e., when
$M_G$ in char 0 specializes to $M_G$ in char $p$. One has to show
that the invariants of the slice representations in char 0
specialize to the invariants in char $p$. However for
$G$=SL(n),the situation is much simpler. One can write down the
automorphism group of a polystable bundle and its representation
on the local moduli space. Consequently, one expects the moduli
spaces to specialize correctly and that the local rings
 of $M_G$ are strongly $F$-regular in all positive characteristics.

We briefly discuss Pic $M_G$ in char 0.  It follows from [4,10]
that $M_G$ has the following properties in char 0:

\begin{enumerate}
\item Pic $M_G \simeq \ZZ$.
\item $M_G$ is a normal projective, Gorenstein variety with rational
singularities and with $K$ negative ample.   \end{enumerate}

Now let $X$ be a normal,Cohen-Macaulay variety in char 0. It is
proved in [13],in response to a conjecture of Karen Smith, that if
$X$ has rational singularities, then the reduction of $X$ mod $p$
is $F$-rational for all large $p$. This result together with 1 and
2 above imply that $M_G$ reduced mod $p$ is $F$-split for all
large $p$. We cannot give effective bounds on the primes involved.
One partial result is known in this direction[12].

 \noindent {\bf Acknowledgement:} I would like to thank my colleagues
 S.~Ilangovan, A.J.~Parame\-swaran and S.~Subramanian for their help
 in preparing this report and T.T.~Nayya and H for constant help and
encouragement.


\label{lastpage}


\begin{thebibliography}{99}

\bibitem{1} M.F. Atiyah, R. Bott, The Yang-Mills equations over Riemann
Surfaces, \emph{Phil. Trans. R. Soc.} London A 308 (1982),
523--615.

\bibitem{2} V. Balaji, A.J. Parameswaran Semistable Principal Bundles-II (in
positive characteristics) to appear in \emph{Transformation
Groups}.

\bibitem{3} V. Balaji, C.S. Seshadri Semistable Principal Bundles-I (in
characteristic zero), to appear in \emph{Journal of Algebra}.

\bibitem{4} A. Beauville, Y. Laszlo, C. Sorger, The Picard Group of the
moduli of $G$-bundles over curves, \emph{Compositio Math.} 112,
(1998), No.2, 183--216.
\bibitem{5} K. Behrend, Semistability of reductive group schemes over curves,
\emph{Math. Ann.} 301 (1995), 281--305.

\bibitem{6} G. Faltings, Stable $G$-bundles and projective connections,
\emph{J. Algebraic Geom.} 2 (1993) No.3, 507--568.

\bibitem{7} D. Gieseker, Stable Vector Bundles and the Frobenius morphism,
\emph{Ann. Sci. Ecol. Nor. Sup.} 6, (1973).

\bibitem{8} M. Hashimoto, Good filtrations of symmetric algebras and strong
$F$-regularity of invariant subrings, \emph{Math. Z.} 236 (2001),
No.3, 605--623.

\bibitem{9} S. Ilangovan, V.B. Mehta, A.J. Parameswaran, Semistability and
Semisimplicity in representations of low-height in positive
characteristics, preprint.

\bibitem{10} S. Kumar, M.S. Narasimhan, Picard group of the moduli spaces of
$G$-bundles, \emph{Math. Ann.} 308, (1997), No.1, 155--173.

\bibitem{11} V.B. Mehta, A.J. Parameswaran, Geometry of low-height
representations, \emph{Proceedings of the International Colloquium
on Algebra, Arithmetic and Geometry}, (ed. R. Parimala), TIFR
Mumbai 2000.

\bibitem{12} V.B. Mehta, T.R. Ramadas, Moduli of vector bundles, Frobenius
splitting and invariant theory, \emph{Ann. of Math.} (2) 144,
(1996), 269--313.

\bibitem{13} V.B. Mehta, V. Srinivas, A characterization of rational
singularities, \emph{Asian J. Math.}, Vol.1, (1997), No.2,
249--271,.

\bibitem{14} V.B. Mehta, S. Subramanian, On the Harder-Narasimhan Filtration
of Principal Bundles, \emph{Proceedings of the International
Colloquium on Algebra, Arithmetic and Geometry}, (ed. R.
Parimala), TIFR Mumbai 2000.

\bibitem{15} V.B. Mehta, S. Subramanian, Moduli of Principal $G$-bundles on
curves in positive characteristic, in preparation.

\bibitem{16} A. Ramanathan, Moduli for principal bundles over algebraic curves
I, \emph{Proc. Indian Acad. Sci. Math. Sci.}, 106, (1996), No.3,
301--328.

\bibitem{17} A. Ramanathan, Moduli for principal bundles over algebraic curves
II, \emph{Proc. Indian Acad. Sci. Math. Sci.}, 106 (1996), No.4,
421--449.

\bibitem{18} A. Ramanathan,  Moduli for principal bundles, in:
\emph{Algebraic Geometry}, Proceedings, Copenhagen 1978, 527--533,
Lecture Notes in Mathematics vol. 732, Springer.

\bibitem{19} S. Ramanan, A. Ramanathan, Some remarks on the instability flag,
\emph{Tohoku Math. Journal} 36, (1984), 269--291.

\bibitem{20} J-P. Serre, \emph{Moursund Lectures}, University of Oregon
Mathematics Department, notes by W.E. Duckworth (1998).
\end{thebibliography}
\end{document}